\documentclass[final]{iopart}   [12pt]
\usepackage{iopams}
\usepackage{color}

\def\be{\begin{equation}}
\def\ee{\end{equation}}
\def\bea{\begin{eqnarray}}
\def\eea{\end{eqnarray}}
\newcommand{\ben}{\begin{equation*}}
\newcommand{\een}{\end{equation*}}
\newcommand{\fr}{\frac}

\newcommand{\z}{\zeta}
\newcommand{\pr}{\partial}
\newcommand{\hs}{\hspace{5mm}}

\newcommand{\dg}{\dagger}

\newcommand{\Ref}[1]{(\ref{#1})}
\newcommand{\PP}{\mathbb P}
\newcommand{\R}{\mathbb R}

\def\chapter#1{\hbox{Introduction.}}

\def\frac#1#2{{\textstyle{#1\over #2}}}

\newcommand{\pa}{\partial}

\begin{document}

\title{Surfaces in ${\mathbb R}^{N^2-1}$ based on harmonic 
maps $S^2\rightarrow CP^{N-1}$}\vskip -0.7cm\hskip 2.5cm\footnote{Talk given at the A. Reznikov memorial workshop at the Max Planck Mathematics Institute in Bonn - 22-29 Sept. 2006}
\vskip 1cm
{\bf Dedicated to the memory of Alexander Reznikov}

\vskip 1cm
\author{W. J. Zakrzewski \footnote{email: w.j.zakrzewski@durham.ac.uk}}
\address{Department of Mathematical Sciences,\\
 University of Durham,
 \\ Durham DH1 3LE, 
 \\United Kingdom}

\date{\today}

\begin{abstract}
\baselineskip=14pt

We show that many surfaces in $\R^{N^2-1}$ can be generated by harmonic maps
of $S^2\rightarrow CP^{N-1}$. These surfaces are based on the projectors 
in $CP^{N-1}$ which describe maps of 
$S^2\rightarrow CP^{N-1}$. In the case when these maps form the Veronese sequence
all the surfaces have constant curvature.

\end{abstract}


\maketitle

\baselineskip=14pt

\section{Introduction}
\label{intro}
\setcounter{equation}{0}
A few years ago, Konopelchenko et al \cite{1,2} initiated a construction
of surfaces immersed in multidimensional 
spaces basing their discussion on the Weierstrass procedure
generalised to higher dimensional spaces.
This has led to further studies \cite{3} and  to relating
the surfaces to the solutions of the  $CP\sp{N-1}$ model\cite{4}.
Recently,  together with Grundland \cite{5}, we have
presented a general procedure for the construction of surfaces
from the harmonic  $CP\sp{N-1}$ maps.
This approach involved writing the equation for the harmonic
map as a conservation law and then observing that the
coordinates of the surfaces can be constructed out of the components
of the special operator which appears in the conservation law.

Our procedure has then been generalised to the supersymmetric case
\cite{6}. At this stage it has become clear that, in the 
holomorphic case, the projector in question is proportional to the fundamental
 projector of the holomorphic map. 
 
 This observation has suggested to us to look at other projectors that
 arise in the description of harmonic maps and to use them to construct
 further surfaces. 
 
 In the next section we recall the general construction of 
 $S^2\rightarrow CP^{N-1}$ harmonic maps and then use them to construct
 surfaces in $\R^{N^2-1}$. 
 In the following section we look, in detail, at the case of $CP^2$
 and show that, in contradistinction to the $CP^1$ case, its surfaces
 do not have to be of constant curvature. We calculate this 
 curvature for two classes of such harmonic maps.
 
 Finally, we look at the Veronese sequence and show that in this case
 all surfaces are of constant curvature.
 We finish the paper with a few comments about the generality of our results.

\section{Harmonic maps $S^2\rightarrow CP\sp{N-1}$\cite{7}}

\setcounter{equation}{0}

\hskip 0.8cm  
The $S^2\rightarrow CP\sp{N-1}$ harmonic maps involve maps 
 into $CP\sp{N-1}$, i.e.
\be
C\ni \zeta=\zeta_1+i\zeta_2\mapsto  z = (z^1,....,z^N)\in 
C\sp{N},
\label{q1}
\ee
where the homogeneous coordinates $z=(z\sp{1},...,z\sp{N})$ 
have the following property
\be
\label{one}
z\sim  z'=\lambda z \quad \hbox{for} \quad \lambda \not= 0.
\ee
Here we have chosen to parametrise the $S^2$ by $\z$ and $\bar{\z}$
the complex variables of the plane obtained from $S^2$ by its
stereographic projection.
Exploiting the projective invariance (\ref{one}) we can require that
\be
z\sp{\dagger}\cdot z\,=\,1
\label{q2}
\ee
holds, where ${\dagger}$ denotes hermitian conjugation. However, 
 we are still left with the gauge (phase) invariance
\be
z \to z'=ze^{i\phi},
\label{q3}
\ee
where $\phi$ is a real-valued function.

 It is easiest to define the $S^2 \rightarrow CP\sp{N-1}$ harmonic maps as stationary 
 points of the Lagrangian  whose density \cite{7} is given by
\be
L\,=\,\frac{1}{4}(D_{\mu}z)\sp{\dagger}\cdot D_{\mu}z,\qquad 
z\sp{\dagger}\cdot z\,=\,1,
\label{q4}
\ee
where the covariant derivatives $D_{\mu}$ act on $z:S^2 \rightarrow CP\sp{N}$
according
to the formula
\be
D_{\mu}z\,=\,\partial_{\mu}z \,-\,(z\sp{\dagger}\cdot \pa_{\mu}z)z.
\label{q5}
\ee
Here the index $\mu=1,2$ denotes the components of $\z_{\mu}$. Note that the 
covariant
derivatives $D_{\mu}z$ transform under the gauge transformation \Ref{q3}
\be
D_{\mu}z \rightarrow D_{\mu}z'=(D_{\mu}z)e^{i\phi},
\label{q6}
\ee
so that the dependence on the phase $\phi$ drops out of the Lagrangian density
\Ref{q4} and so the target space is really $CP\sp{N-1}$.

The total Lagrangian is given by
\be
{\cal L}\,=\,\int L\,d\zeta_1d\zeta_2
\label{q7}
\ee
and, as we want to consider $S^2 \rightarrow CP\sp{N-1}$ maps, we require
that ${\cal L}$ is finite.

To find the maps it is convenient to define 
\be
z\,=\,{f\over \vert f\vert},
\label{q8}
\ee
where $|f|=(f\sp{\dg}\cdot f)^{1/2}$.
In terms of $f$ the Lagrangian \Ref{q7} becomes
\be
{\cal L}\,=\,\int {|\bar{\pa} f|^2+|\pa f|^2\over |f|^4}\,d\zeta\,d\bar{\zeta}
\label{q9}
\ee
where $|\pa f|^2=(\pa f)\sp{\dg}\cdot (\pa f)$ and
$\vert \bar\pa f\vert\sp2 =(\bar{\pa} f)\sp{\dg} \cdot (\bar{\pa} f)$.
The Euler-Lagrange equations for $f$ take the form
\be
\left(1\,-\,{f\otimes f\sp{\dg}\over \vert f\vert\sp2}\right)\left[\pa\bar\pa f\,-\,\pa f\,{(f\sp{\dg}\cdot \bar\pa f)\over \vert f\vert\sp2}\, -\,\bar\pa f\,{(f\sp{\dg}\cdot \pa f)\over \vert f\vert\sp2}\right]\,=\,0,
\label{q10}
\ee
where we have introduced the holomorphic and antiholomorphic 
derivatives
\be
\partial ={\partial\over \partial(\zeta_1 +i\zeta_2 )}\,=\,{\pa\over \pa \z},
\label{partia}\quad 
\bar\partial ={\partial\over \partial(\zeta_1 -i\zeta_2)}\,=\,{\pa\over \pa \bar\z}
\ee
and bar denotes complex conjugation.

Then all harmonic maps of $S^2\rightarrow CP\sp{N-1}$ can be constructed in the following
way:

First observe that $f=f(x_+)$, {\it i.e} whose components are analytical functions of 
$\zeta=\zeta_1 +i\zeta_2$ automatically satisfies (\ref{q10}). For the total Lagrangian ${\cal L}$
to be finite we require that each component of $f$ is a polynomial in $\zeta$. Then any such 
$f$ describes a harmonic map of $S^2\rightarrow CP^{N-1}$ and as $f$ is analytic
such a map is also holomorphic.

Other maps can be constructed from the holomorphic maps. 
To do this we define an operator $P_+$  by its action in $C^{N}$
\be
P_+\, g\,= \partial g\, -\, g\,{g^{\dagger}\cdot \partial g\over g^{\dagger}\cdot g}.
\ee

Then we apply this operator to our holomorphic vector $f$ obtaining $P_+f$. And repeat 
this procedure applying it to $P_+f$ obtaining $P_+^2f$ and so on.
Then, as is known \cite{7},  all $P_+^kf$ solve (\ref{q10}) and so represent
harmonic maps $S^2\rightarrow CP^{N-1}$.

The sequence of vectors $P_+^kf$ have the following properties \cite{7}:
\
\begin{eqnarray} 
(P_+^if)^{\dagger}\cdot P_+^jf\,=\,0,&&\quad i\ne j, \nonumber\\
\bar\pr\left(P^k_+ f\right)=-P^{k-1}_+ f \fr{|P^k_+
f|^2}{|P^{k-1}_+ f|^2},
\hs &&
\pr\left(\fr{P^{k-1}_+ f}{|P^{k-1}_+ f|^2}\right)=\fr{P^k_+
f}{|P^{k-1}_+f|^2}.
\label{prop}
\end{eqnarray} 

Next we construct a sequence of  projectors  $P_k$ 
of the form 
\be
P(V)={V \otimes V^\dg \over |V|^2},
\label{for}
\ee
where for $V$ we take
$V=P\sp{k}_+f$.

Due to the properties (2.14) these projectors are mutually orthogonal
and we have
\be
\sum_{k=0}^{N-1}\,P_k\,=\,1.
\ee
 Hence only $N-1$ of $P_k$ are independent ({\it ie} we can consider
only those corresponding to $k=0,1,..,N-2)$.

\section{Surfaces in $\R^{N^2-1}$}

\setcounter{equation}{0}
 
The Weierstrass construction of surfaces discussed in \cite{1} involves looking at a different
set of equations and then using their solutions to construct surfaces.
However, this set of equations was shown in \cite{4} to be equivalent 
to the harmonic maps so one could use these maps directly - to construct our surfaces.
In fact in \cite{5} it was shown how to construct such surfaces and, at the same time,
the generalised Weierstrass system was given.

The work in \cite{6} has revealed that, in the $CP^1$ case, the surfaces are 
related to the 
basic projector of the corresponding map. This suggests that we look at other 
projectors and use them to construct further surfaces.

The orthogonality of the projectors allows us to take their linear combinations.
Hence, in the $CP^{N-1}$ case we can take as our projector
\be
\PP\,=\,\sum_{k=0}^{N-2}\,\alpha_k\,P_k,
\ee
where $\alpha_k$ are constants.
From this matrix $\PP$ we construct a vector $X$, with $N^2-1$ components, in the 
following way:

First we consider the off-diagonal entries of $\PP$ (note that $\PP$ is hermitian) and use them to define $N^2-N$ 
components of $X$ by taking their real and imaginary parts; {\it ie} we take
\be
X_l\,=\,\PP_{ij}+\PP_{ji},\quad {\hbox{for}}\quad l=1,...,\frac{N^2-N}{2},\quad i\ne j
\ee
$$
X_l\,=\,i(\PP_{ij}-\PP_{ji}),\quad {\hbox{for}}\quad l=\frac{N^2-N}{2},...,{N^2-N},\quad i\ne j.
$$
The remaining $N-1$ components of $X$ are taken from the diagonal entries of $\PP$.

As each projector $P_k$ has trace one, the trace of $\PP$ is also constant so that the different
choices of these $N-1$ components of $X$ correspond to the shifts of the vector $X$ and so would not
alter the metric nor the curvature of the surface in $\R^{N^2-1}$.

How do we choose these remaining components of $X$?
The obvious procedure is to choose them in such a way that
\be
\sum_{i=1}^{N-1}\partial_+X_i\partial_-X_i\,=\, 2\, \sum_{i=0}^N\partial_+{\mathbb P}_{ii}\partial_-{\mathbb P}_{ii}.
\ee
In the ${\mathbb C}P^1$ case this tells us that 
for the last component of $X$ we should take 
 ${\mathbb P}_{11}-{\mathbb P}_{22}$. For larger $N$ we have more choices; thus,
as was discussed in \cite{6}, for  $CP\sp{2}$  we can 
take (this choice is based on Gell Mann's $SU(3)$ $\lambda$ matrices)
\be 
X_1\,=\,{\mathbb P}_{11}\,-\,{\mathbb P}_{22}, \qquad  X_2\,=\,\sqrt{3}({\mathbb P}_{11}\,+\,{\mathbb P}_{22}).
\ee
or we could make another choice. In general, for  $CP\sp{2}$,  we could take
\be
{\mathbb P}_{11}\,=\, \frac{1}{3}\,+\,a X_1\,+\, bX_2,\qquad {\mathbb P}_{22}\,=\,
 \frac{1}{3}\,+\,cX_1\,+\, dX_2.
\ee

Then we choose $a$, $b$, $c$ and $d$ so that
\be 
\partial_+X_{1}\partial_-X_{1}+\partial_+X_{2}\partial_-X_{2}
\ee
gives the same expression as
\be
\partial_+{\mathbb P}_{11}\partial_-{\mathbb P}_{11}+\partial_+{\mathbb P}_{22}\partial_-{\mathbb P}_{22}+
\partial_+{\mathbb P}_{33}\partial_-{\mathbb P}_{33}
\ee
in which we can eliminate ${\mathbb P}_{33}$ using ${\mathbb P}_{33}=1-{\mathbb P}_{11}-{\mathbb P}_{22}.$

A simple calculation shows that we have a one-parameter family of solutions
\be
a\,=\,\frac{2}{\sqrt{3}}\cos\alpha,\quad b\,=\,\frac{2}{\sqrt{3}}\sin\alpha,
\ee

$$
c\,=\,\mp\sin\alpha\,-\frac{1}{\sqrt{3}}\cos\alpha\,,\quad d\,=\,-\frac{1}{\sqrt{3}}\sin\alpha\,\pm \cos\alpha.
$$

For $N>2$ the solutions are even more nonunique.
 
\section{Properties of the surfaces}
\label{prop}
\setcounter{equation}{0}

Let us consider a surface defined by $X$ (this is clearly a surface as our vector  $X$ depends 
on two variables $\z$ and $\bar{\z}$).

The metric on the surface, induced by the map, is, due to our choice of $X$, given by
\be
g_{++}\,=\,\tr \,({\pa \PP \pa \PP}),\quad g_{+-}\,=\,\tr \,({\pa \PP \bar\pa \PP})
\ee
and, of course, $g_{--}=\overline{g_{++}}$.

To calculate the metric we need some properties of the projectors $P_k$.

Note that due to (\ref{prop}) we have
\be
\pa P_k \,=\,{P_+^{k+1}f \otimes (P_+^kf)^\dg \over |P_+^kf|^2} \,-\, {P_+^kf \otimes (P_+^{k-1})^\dg \over |P_+f^{k-1}|^2}
\ee
(for $k=0$ we have only the first term)
and
$\bar {\pa} P_k\,=\,\overline{\pa P_k}$.

Then, due to the orthogonalisty of $P_+^kf$ we see that
\be 
g_{++}\,=\,g_{--}\,=\,0
\ee
and that
\be
g_{+-}\,=\, \alpha_0{|P_+f|^2\over |f|^2}\,+\,\sum_{k=1}^{N-2}\,\alpha_k\left({|P_+^{k+1}f|^2 \over |P_+^{k}f|^2}
+{|P_+^{k}f|^2\over| P_+^{k-1}f|^2}\right).
\ee

As only the $g_{+-}$ component of the metric is nonzero the curvature is given by \cite{8}
\be
K\,=\,-{4\over g_{+-}}\,\pa\bar{\pa}\,\ln(g_{+-}).
\ee
In general, it is difficult to calculate the curvature for the expression above; hence in the next
section we discuss its form in some special cases.

\section{Special cases}
\label{prop}
\setcounter{equation}{0}
\subsection{$CP^1$}

Consider first the case of $CP^1$ maps. In this case have only holomorphic harmonic maps
and the surface in $\R^3$ is a sphere \cite{9}. To see this we note that we can take 
\be
X_1\,=\,{W+\bar W\over 1+|W|^2},\quad X_2\,=\,i{W-\bar W\over 1+|W|^2},\quad X_3\,=\,{2\over 1+|W|^2},
\ee
and this gives us
\be
X_1^2+X_2^2+(X_3-\frac{1}{2})^2\,=\,\frac{1}{4}.
\ee

The induced metric is 
\be 
g_{+-}\,=\,{|P_+f|^2\over |f|^2}.
\ee
Note that as $f$ has only two components we can put $f=(1,W)$ where $W$ is a ratio of polynomials 
in $\z$. Then 
\be g_{+-}={|\pa W|^2\over (1+|W|^2)^2}
\ee  and
\be
K\,=\,-4{(1+|W|^2)^2\over |\pa W|^2} \pa \bar{\pa} \, \ln({|\pa W|^2\over (1+|W|^2)^2})
\,\ee
$$=\,4 {(1+|W|^2)^2\over |\pa W|^2} \pa \bar{\pa} \, \ln({ (1+|W|^2)^2})\,=\,8 {(1+|W|^2)^2\over |\pa W|^2}
{|\pa W|^2\over (1+|W|^2)^2}\,=\,8 $$
as $\bar{\pa} W=0$.
This is, of course, the curvature of the $CP^1$ harmonic map and is, also, the surface
generated by Konopelchenko via his Weierstrass procedure.

\subsection{$CP^2$ case}

Now we have more choices. We have two classes of harmonic maps; the holomorphic ones
({\it ie} those based on $f$) and nonholomorphic ones (based on $P_+f$).
Our vector $X$ has 8 components and its entries are constructed from the matrix $\PP$ which in this case
takes the form
\be 
\PP\,=\,\alpha_o P_0 \,+\,\alpha_1 P_1.
\ee
(of course, in addition we has some freedom of how to choose the two compenents of $X$
constructed out of the diagonal entries of $\PP$).
The metric is now given by
\be
g_{+-}\,=\,(\alpha_0+\alpha_1){|P_+f|^2\over |f|^2}\,+\,\alpha_1{|P_+^2f|^2\over |P_+f|^2}.
\ee
Let us look first at the case of the holomorphic map ({\it ie} when $\alpha_1=0$ and $\alpha_0=1$).
Then, $g_{+-}={|P_+f|^2\over |f|^2}$ and the curvature is given by
\be
K\,=\,4\left(2\,-\,{|P_+^2f|^2\,|f|^2\over |P_+f|^4}\right).
\label{aa}
\ee
Note that this curvature is, in general, not constant.
In the cases of the embedding of $CP^1$ into $CP^2$  $f$ has only 2 components and then 
$P_+^2f=0$ and so the result reduces back to the $CP^1$ case.

The curvature is also constant when the second term in (\ref{aa}) is constant and this case
corresponds to the Veronese sequence. We shall discuss this case later.

Next we look at the nonholomorphic case, {\it ie} when $\PP=P_1$.
Then 
\be
g_{+-}\,=\,{|P_+f|^2\over |f|^2}\,+\,{|P_+^2f|^2\over |P_+f|^2}
\ee
and the calculation of the curvature is quite tedious. In fact, lengthy calculations give
\be
\label{cp2}
K\,=\,4{2\left({|P_+f|^6\over |f|^6}\,+\,{|P_+^2f|^6\over |P_+f|^6}\right)
\,+\,{|P_+^2f|^2\over |f|^2}\left|\pa\,\ln\,{|f|^2|P_+^2f|^2\over |P_+f|^4}\right|
\over \left({|P_+f|^2\over |f|^2}\,+\,{|P_+^2f|^2\over |P_+f|^2}\right)^3}.
\ee
Of course if we take a more general case (with both $\alpha_0$ and $\alpha_1$ nonvanishing)
we get an even more complicated expression.

\subsection{Veronese sequence}

The Veronese choice of the vector $f$ is such that its components are monomials of $\z$ multiplied by the square roots of 
the coefficients of the expansion of $(1+a)^N$ in powers of $a$.
Hence $f$ is given by
\be
f=(1,\,\sqrt{N}\z,\, \sqrt{\frac{N(N-1)}{2}}\z^2,\, ... \sqrt{N}\z^{N-2},\,\z^{N-1})
\ee
With this choice all $|P_+^kf|^2$ are given by the powers of $(1+|\z|^2)$.
Moreover, $|f|^2=(1+|\z|^2)^{N-1}$ and the successive powers decrease by 2.
Thus
\be
g_{+-}\,=\,\alpha\,\left({1\over 1+ |\z|^2}\right)^2
\ee
for all harmonic maps. Only the value of $\alpha$ depends on the map.

Hence the curvature is given by
\be
K\,=\,{8\over \alpha}.
\ee

Looking at the concrete examples we see that the surfaces are quite complicated.
For example, if we consider the $CP^2$ case and look at the two projectors $P_0$ and $P_1$ we note 
that the corresponding surfaces are very different.

This is clear as $P_0$ is generated by $f$ which in this case is given by $f=(1,\sqrt{2}\z,\z^2)$ and so
\be
\label{P0}
P_0\,=\,{1\over (1+|\z|^2)^4}\left(\begin{array}{ccc}
1&\sqrt{2}\bar\z&\bar\z^2\\
\sqrt{2}\z&2|\z|^2&\sqrt{2}\bar\z|\z|^2\\
\z^2&\sqrt{2}\z|\z|^2&|\z|^4\end{array}\right).
\ee
On the other hand \be
P_+f\,=\,{\sqrt{2}\over 1+|\z|^2}\left(-\sqrt{2}\bar\z,\,(1-|\z|^2,\,\sqrt{2}\z\right)
\ee
and so
\be
\label{P1}\!\!\!\!\!\!\!\!\!\!\!\!\!\!\!\!
\!\!\!\!\!\!\!\!\!\!\!P_1\,=\,{1\over (1+|\z|^2)^2}\left(\begin{array}{ccc}
2|\z|^2&-\sqrt{2}\bar\z(1-|\z|^2)&-2\bar\z^2\\
-\sqrt{2}\z(1-|\z|^2)&(1-|\z|^2)^2&\sqrt{2}\bar\z(1-|\z|^2)\\
-2\z^2&\sqrt{2}\z(1-|\z|^2)&2|\z|^2\end{array}\right).
\ee

Note that as several components of $P_1$ are proportional to each other the corresponding vector $X$ has
some components the same and so changing the basis in $\R^8$ we observe that vector $X$ lies in a 
five-dimensional subspace of $\R^8$ and so we can take it in the form
$$
X_1\,=\, {2x(1-x^2-y^2)\over (1+x^2+y^2)^2},\,\quad X_2\,=\, {2y(1-x^2-y^2)\over (1+x^2+y^2)^2},
$$
\be
X_3\,=\,{2(x^2-y^2)\over (1+x^2+y^2)^2},\,\qquad X_4\,=\,{4xy\over (1+x^2+y^2)^2}
\ee
$$X_5\,=\,\sqrt{3}{(1-x^2-y^2)^2\over (1+x^2+y^2)^2},
$$
where we have defined $x$ and $y$ through $\z=x+iy$.

Note that the curvatures for the two cases ($P_0$ and $P_1$) are constant but different; namely:
\be
K(P_0)\,=\, 4,\qquad K(P_1)\,=\,2.
\ee

\section{Conclusions}

We have discussed here a possible construction of surfaces in $\R^{N^2-1}$ based on harmonic maps
$S^2\rightarrow CP^{N-1}$. Our construction, which in a way, is a generalisation of the Weierstrass
construction used by Konopelchenko and collaborators has produced many surfaces whose induced metric
is related to the total Lagrangian (energy of the underlying maps).
The curvatures of these surfaces can be easily calculated. In the $\R^3$ case - the surface is a sphere but 
for larger $N$ the surfaces are more complicated. Moreever, their curvatures are related to the curvatures
of the $CP^{N-1}$ spaces. When we restrict our attention to the Veronese sequence of maps
all corresponding surfaces have a constant curvature.

\label{concl}
\setcounter{equation}{0}

\section*{Acknowledgments}
The work reported in this paper was prepared for the Max Planck Symposium 
dedicated to the memory of Alexander Reznikov, who was my colleague in Durham
and with who I discussed the topics of this paper. I am grateful to the 
organisors for  inviting me to the symposium and giving me an opportunity
to present this paper.

The results reported in this paper are the natural outgrowth
of the work done with V. Hussin, A.M. Grundland and A. Strasburger. I would like to thank them
for their collaboration. I also would like to thank J. Bolton for 
his interest and for helpful
discussions.


\section*{References}

\begin{thebibliography}{20}
\bibitem{1}
B. Konopelchenko and I. Taimanov, Constant mean curvature surfaces via
an integrable dynamical system, {\it J. Phys.} {\bf A 29}, 1261-1265 (1996).

\bibitem{2}
R. Carroll and B. Konopelchenko, Generalised Weierstrass-Enneper inducing
conformal immersions and gravity, {\it Int. J. Mod. Phys.} {\bf A 11}, (7), 
1183-1216 (1996).

\bibitem{3}
B. Konopelchenko and G. Landolfi, Generalised Weierstrass representation
for surfaces in multi-dimensional Riemanian spaces, {\it Stud. Appl. Maths.}
{\bf 104}, 129-169 (1999) and references therein.

\bibitem{4}
 P. Bracken and A.M. Grundland, Symmetry properties and explicit solutions
of the generalised Weierstrass system, {\it J. Math. Phys.} {\bf 42}, 1250-1282
(2001) and references therein.

\bibitem{5}
A.M. Grundland, and W.J. Zakrzewski, 
$CP^{N-1}$ harmonic maps and the Weierstrass
problem, {\it J. Math. Phys.} {\bf 44}, 3370-3382 (2003).

\bibitem{6}
V. Hussin and W.J. Zakrzewski, Susy $CP\sp{N-1}$ model and surfaces in ${\mathbb R}^{N^2-1}$,  preprint to be published in {\it J. Phys.} {\bf A} (2006). 

\bibitem{7}
 see {\it eg} W.J. Zakrzewski, {\em Low Dimensional Sigma Models} (Hilger, Bristol, 1989).

\bibitem{8} J. Bolton, G.R. Jensen, M. Rigoli and L.M. Woodward, 
On conformal minimal immersions of $S^2$ into $CP^{n}$, {\it Math. Ann.} 
{\bf 279}, 599-620 (1988)

\bibitem{9} E.V. Ferapontov and A.M. Grundland, Links between different analytic descriptions of constant mean curvature surfaces, {\it J. Nonlin. Math. Phys.}
{\bf 7}, 14-21 (2000)
\end{thebibliography}
\end{document}